\documentclass[12pt]{article}
\usepackage{calc,graphics}
\usepackage{multirow}
\usepackage{todonotes}
\usepackage{booktabs}
\usepackage{pdfpages}
\usepackage{times}
\usepackage{bm}
\usepackage{comment}
\usepackage{amsgen,amsmath,amssymb,latexsym}
\usepackage{graphicx}
\usepackage{enumerate}
\usepackage{type1cm}
\usepackage{xspace}
\usepackage{tikz}
\usepackage{color}
\usepackage{rotating}
\usetikzlibrary{shapes,arrows,shadows,decorations}
\usetikzlibrary{decorations.pathmorphing}
\usetikzlibrary{backgrounds}

\newcommand{\newpar}{{\vspace{0.3cm} \noindent}}
\newcommand{\indep}{\mbox{$\,\perp\!\!\!\perp\,$}}
\newcommand{\nindep}{\mbox{$\,\not\!\perp\!\!\!\perp\,$}}

\newcommand{\etal}{{\em et al\/}\xspace}

\newtheorem{condition}{Condition}

\DeclareGraphicsExtensions{.pdf,.png,.gif,.jpg,.tif,.ps}

\begin{document}

\title{\huge
Bayesian Mendelian Randomization}
\author{
\rule{0cm}{1cm}
Carlo Berzuini$^{1}$,
Hui Guo$^{1}$, Stephen Burgess$^{2}$,
Luisa Bernardinelli$^{3}$\\
\rule{0cm}{0.6cm}\\
\begin{minipage}{15cm}
{\large
\noindent $^{1}$ Centre
for Biostatistics, The University of Manchester,
Manchester,
UK\thanks{carlo.berzuini@manchester.ac.uk}\\
\noindent $^{2}$ Department of Public
Health and Primary Care, University of Cambridge, Cambridge, United Kingdom\\
\noindent $^{3}$ Department of Brain and
Behavioural Sciences, University of Pavia,
Italy
}
\end{minipage}
}
\date{}
\maketitle

\newpar {\bf Keywords:} Pleiotropy;
Horseshoe Prior;
Shrinkage Estimators; 
Egger Regression; Median
Estimator; 
Causal Inference; Instrumental Variable Analysis;
Invalid Instruments;
Instrument-Exposure Interaction; Metabolomics


\begin{abstract}
\noindent Our Bayesian approach to Mendelian Randomization
uses multiple instruments to assess the putative causal effect of
an exposure on an outcome. The approach is robust to violations
of the (untestable) Exclusion Restriction condition, and hence it does not require
instruments to be independent of the outcome conditional on the
exposure and on the confounders of the exposure-outcome relationship.
The Bayesian approach offers a rigorous handling of
the uncertainty (e.g. about the estimated instrument-exposure
associations), freedom from asymptotic
approximations of the null distribution
and the possibility to elaborate the model
in any direction of scientific relevance.
We illustrate the last feature with the aid
of a study of the metabolic mediators of the
disease-inducing effects of obesity,
where we elaborate the model to investigate
whether the causal effect of interest interacts with
a covariate. The proposed model contains a
vector of unidentifiable
parameters, $\beta$, whose $j$th
element represents the pleiotropic
(i.e., not mediated by the exposure)
component of the association of instrument $j$
with the outcome.  We deal with
the incomplete identifiability
by assuming that the pleiotropic effect
of some instruments
is null, or nearly so, formally by
imposing on $\beta$
Carvalho's horseshoe shrinkage prior,
in such a way that different
components of $\beta$
are subjected to different degrees
of shrinking, adaptively and in
accord with the compatibility
of each individual instrument
with the hypothesis of no pleiotropy.
This prior requires a minimal input from
the user. We present the results of a simulation
study into
the performance of the proposed method
under different types of
pleiotropy and sample sizes. Comparisons
with the performance of the weighted median
estimator are made. Choice of the prior and
inference via Markov chain Monte Carlo
are discussed.
\end{abstract}

\section{Introduction}
\label{Introduction}

Mendelian randomization (MR) is a method for testing and estimating the causal effect
of an exposure, $X$, on an outcome, $Y$ in situations where the relationship between
these two variables is confounded, based on information from the genotyping of variants
associated with the exposure \cite{Burgess2012} \cite{Davey2003}  \cite{Davey2014}. The
method assumes that the genotypes are the result of a randomised experiment performed
by Nature during meiosis, and therefore their effect on the outcome is likely to mimick
the effect we would obtain through an intervention on the exposure. Standard MR
theory requires the instruments to be  {\em (i)} associated with the exposure, {\em (ii)} 
independent of the outcome conditional on the exposure and on the confounders of the
exposure-outcome relationship, and {\em (iii)} independent of those confounders.
Many authors, e.g. Jones, Didelez and colleagues \cite{Didelez2007} 
\cite{Didelez2007b} \cite{Jones2012} call {\em (ii)} the Exclusion Restriction condition.

\newpar Current, state-of-the-art, MR methods combine the strengths of multiple instrumental
variants, $Z =(Z_1, \ldots , Z_J)$, and are robust to violations of the Exclusion
Restriction condition {\em (ii)}. Examples include
 the Egger regression and the median estimator method \cite{Bowden2016} \cite{burgess2016sensitivity}
\cite{Burgess2015}, both of which work from the estimated coefficient, $\beta_{YZ_j}$, of the
regression of $Y$ on each $Z_j$, and the estimated coefficient, $\beta_{XZ_j}$, of the regression
of $X$  on each $Z_j$. The Egger regression method interprets a linear relationship between the
$\{\beta_{YZ_j}\}$ and the $\{\beta_{XZ_j}\}$ as evidence that an exogenous change in $X$  will
cause a corresponding proportional change in $Y$, and allows any proportion of the $J$ instruments
to violate the Exclusion Restriction condition \cite{Bowden2015}, whereas the median
estimator assumes that
only a minority of the instruments violate that condition.

\newpar One limitation of the median and Egger estimators is that they both treat the
estimated coefficient of the regression of $X$ on each instrument as fixed, despite the
considerable uncertainty that may surround it. This confronts the user with a problematic
trade-off between using all the available data information and excluding the weakly associated
variants.

\newpar This paper presents a Bayesian approach to MR that uses multiple
instruments and relaxes the Exclusion Restriction condition. Not only does this
approach overcome the above mentioned limitations; it also offers
the typical advantages of Bayesian analysis, including a rigorous
handling of the uncertainty, the freedom from asymptotic
approximations of the null distribution and the possibility
of a straightforward elaboration of the model in any direction of scientific relevance.
We illustrate the last feature  in Section \ref{Sex-Dependent Causal Effect of
Body Mass on Phenylalanine}, with the aid of a study of the metabolic mediators of the
disease-inducing effects of obesity, where we elaborate the model to investigate
whether the causal effect of obesity on certain mediators interacts with sex.

\newpar We start our journey in Section \ref{Causal Assumptions} by
specifying the causal assumptions behind the
method. An additional assumption is introduced in
Section \ref{The Likelihood}. The model proposed in Section
\ref{The Model} contains a vector of unidentifiable parameters,
$\beta$,  whose $j$th element represents the pleiotropic
(i.e., not mediated by the exposure) component of the association
of instrument $j$ with the outcome. In the same Section we show
that the incomplete identifiability can be dealt with by assuming
that the pleiotropic effect of some instruments is null, or nearly so. 
We do this by imposing on $\beta$ the horseshoe shrinkage prior proposed
by Carvalho and colleagues \cite{carvalho2010}, under which
different components of $\beta$ will be subjected to different degrees
of shrinkage, inferred from the data. The horseshoe
prior requires a minimal input from the user. This and other prior
specifications are discussed in Section \ref{The Prior}.

\newpar In Section \ref{Simulation Experiment} we perform
a simulation experiment to assess
the performance of the method in the presence of
different types of pleiotropy and with different sample sizes, 
and to compare it with that of the weighted median estimator. 
The inference is performed by sampling
the model posterior distribution via the
Hamiltonian dynamics Markov chain
Monte Carlo techniques \cite{Metropolis1953} \cite{Neal2011}
incorporated in the program {\tt Stan} \cite{STAN2014b}
\cite{STAN2014a}.
Sensible initial values for the chains
are obtained from estimates of the
posterior means
obtained via
variational inference techniques 
\cite{Wainwright2008}.

\newpar This work restricts attention to
the case where $X$ and $Y$ are continuous. The discrete
case is under investigation.

\begin{center}
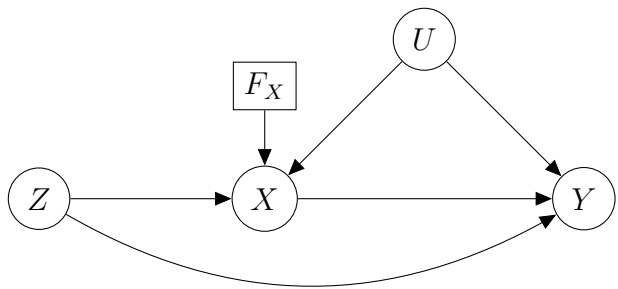
\begin{figure}[!h]
\centering
\begin{tikzpicture}[align=center,node distance=1.5cm]
\node (X) [draw=black,fill=none,circle] {$X$};
\node (U) [above right of=X,draw=black,fill=none,circle,
node distance=3cm] {$U$};
\node (Y) [below right of=U,draw=black,fill=none,
circle,node distance=3cm] {$Y$};
\node (Z) [left of=X,draw=black,fill=none,circle,
node distance=3cm] {$Z$};
\node (SX) [above of=X,draw=black,fill=none,rectangle]{$F_X$};
\draw[-triangle 45] (X) to node
[auto] {} (Y);
\draw[-triangle 45] (U) to node {} (X);
\draw[-triangle 45] (U) to node {} (Y);
\draw[-triangle 45] (Z) to node {} (X);
\draw[-triangle 45] (SX) to node {} (X);
\draw[-triangle 45, bend right] (Z) to node {} (Y);
\end{tikzpicture}\label{Figure 1}
\caption{\small Graphical representation of the assumptions
at the basis of the proposed method.}
\end{figure}
\end{center}

\begin{center}
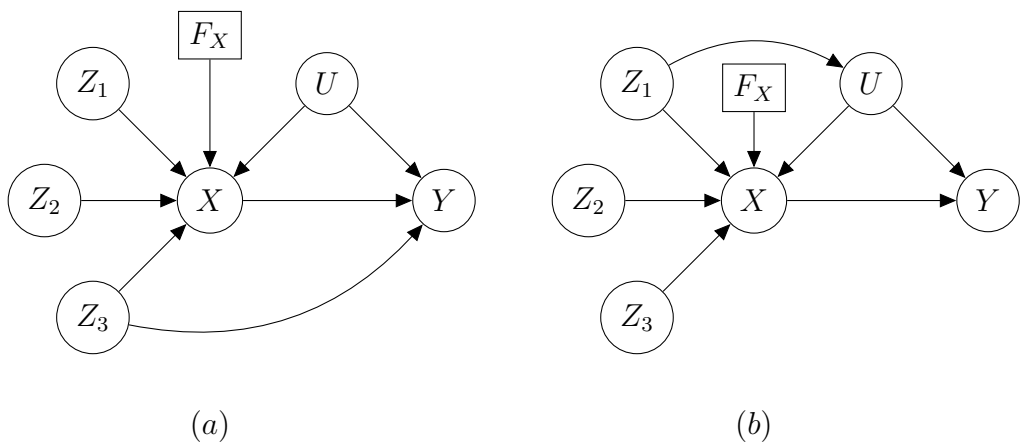
\begin{figure}[!h]
\centering
\begin{tikzpicture}[align=center,node distance=2.2cm]
\node (X) [draw=black,fill=none,circle] {$X$};
\node (U) [above right of=X,draw=black,fill=none,circle] {$U$};
\node (Y) [below right of=U,draw=black,fill=none,
circle] {$Y$};
\node (Z1) [above left of=X,draw=black,fill=none,circle] {$Z_1$};
\node (Z2) [left of=X,draw=black,fill=none,circle] {$Z_2$};
\node (Z3) [below left of=X,draw=black,fill=none,
circle] {$Z_3$};
\node (a) [below of=X, node distance=3cm] {$(a)$};
\node (SX) [above of=X,draw=black,fill=none,rectangle]{$F_X$};
\draw[-triangle 45] (X) to node
[auto] {} (Y);
\draw[-triangle 45] (U) to node {} (X);
\draw[-triangle 45] (U) to node {} (Y);
\draw[-triangle 45] (Z1) to node {} (X);
\draw[-triangle 45] (Z2) to node {} (X);
\draw[-triangle 45] (Z3) to node {} (X);
\%draw[-triangle 45] (S3) to node {} (Z3);
\draw[-triangle 45] (SX) to node {} (X);
\draw[-triangle 45, bend right] (Z3) to node {} (Y);
\end{tikzpicture} \hspace{0.8cm} \begin{tikzpicture}[align=center,node distance=2.2cm]
\node (X) [draw=black,fill=none,circle] {$X$};
\node (U) [above right of=X,draw=black,fill=none,circle] {$U$};
\node (Y) [below right of=U,draw=black,fill=none,
circle] {$Y$};
\node (Z1) [above left of=X,draw=black,fill=none,circle] {$Z_1$};
\node (Z2) [left of=X,draw=black,fill=none,circle] {$Z_2$};
\node (Z3) [below left of=X,draw=black,fill=none,circle] {$Z_3$};
\node (b) [below of=X,node distance=3 cm] {$( b)$};
\node (SX) [above of=X,draw=black,fill=none,rectangle,node distance=1.5cm]{$F_X$};
\draw[-triangle 45] (X) to node [auto] {} (Y);
\draw[-triangle 45] (U) to node {} (X);
\draw[-triangle 45] (U) to node {} (Y);
\draw[-triangle 45] (Z1) to node {} (X);
\draw[-triangle 45] (Z2) to node {} (X);
\draw[-triangle 45] (Z3) to node {} (X);
\draw[-triangle 45] (SX) to node {} (X);
\draw[-triangle 45,bend left] (Z1) to node {} (U);
\end{tikzpicture}
\label{Figure 3}
\caption{\small The three instruments in example {\em (a)} satisfy Condition \ref{confounder independence},
whereas instrument $Z_3$ in example {\em (b)} violates it.}
\end{figure}
\end{center}

\section{Causal Assumptions}
\label{Causal Assumptions}

\newpar Let the symbol $Z$, with $Z \equiv (Z_1, \ldots , Z_J)$, represent a set of instruments,
and $U$ the set of (generally unknown) common causal influences of $X$ and  $Y$.
Let the notation $A \indep B \mid C$ hereafter stand for "$A$ independent of $B$
given $C$" \cite{Dawid1979}.  Within our method,
 $Z$ qualifies as a set of instruments for assessing the putative causal effect of $X$ on $Y$ 
if there exists a (possibly empty) set of observed variables $W$ such that, conditional on $W$,
each instrument $Z_j$ satisfies the following conditions:

\begin{condition}[relevance]
\label{relevance}
the instrument is associated (not necessarily in a causal way) with
the exposure:
$Z_j \nindep X$.
\end{condition}
\begin{condition}[confounder independence]
\label{confounder independence}
the instrument is independent of the confounders $U$: $Z_j \indep U$.
\end{condition}
\noindent We say that the variables
in $Z$  {\em unconditionally} qualify as instruments within our method
if the above conditions hold with $W$ empty. Conditioning on $W$ will 
hereafter be taken as implicit in the notation.  Condition
\ref{confounder independence} is untestable, as it involves the unknown quantity $U$.
Our method is robust to violations of the (untestable) Exclusion Restriction condition,
 $Y \indep Z \mid (X,U)$, insofar as $Z \indep U$ remains valid. In
particular, the instruments are allowed to influence $Y$ through
pathways independent of $U$ and $X$. This is an important generalisation, if
one considers how difficult it is to corroborate the Exclusion Restriction
condition, whether on the basis of empirical evidence or
biological knowledge.

\newpar A general situation where the set of instruments $Z$ 
unconditionally satisfies Condition \ref{confounder independence}
is depicted by the conditional independence graph of Figure 1 where,
for the time being, the reader is asked to ignore node $F_X$. One
can use the $d$-separation rule \cite{geiger90}, or moralisation
 \cite{lauritzen90}, to check that this graph
violates the Exclusion Restriction condition, due to the presence
of a $Z \rightarrow Y$ arrow,
but this condition is not required by our method. Hence, if the assumptions
in Figure 1 hold and $Z$ is associated with $X$, then $Z$
unconditionally qualifies as an instrument for assessing the
causal effect of $X$ on $Y$ through our method.
While the example of Figure 2{\em (a)} is a special case of
Figure 1, Figure 2{\em (b)} contains an instrument ($Z_3$) that
violates Confounder Independence, which is not
compatible with our method.

\newpar The node $F_{X}$ in Figure 1 represents an example of regime
indicator \cite{Dawid2002} \cite{Dawid2015}, and tells us whether the
value of $X$ is set by a hypothetical exogenous
intervention or instead it arises from passive observation. 
Embodied in the graph is the relationship
$F_X \indep (U,Z)$, stating that an intervention
on $X$ will not affect $Z$ or $U$ -- a sensible
assumption if we accept that genetic variants cannot be causally affected by
changes in $X$. Also expressed in the graph is the
assumption $Y \indep F_X \mid (X,Z,U)$, stating that, conditional
on $Z$ and $U$, the distribution of $Y$ given $X$ does not depend on whether
the value of $X$ has been generated by passive observation
or intervention. This implies that the $X \rightarrow Y$ arrow in the graph,
and the coefficient of $X$ in a regression of $Y$ on $(X,Z,U)$, can be
interpreted causally. By contrast, the method does not necessarily require
the association between $Z$ and $X$ (represented in the graph by the
$Z \rightarrow X$ arrow) to be causal.

\section{The Model}
\label{The Model}

\subsection{The Likelihood}
\label{The Likelihood}

\noindent Assume each sample individual
is characterized by a complete set
of observed values for $X, Y$
and $Z \equiv  (Z_{1},
\ldots , Z_{J})$.  According to the graph in Figure 1,
the conditional distribution
$P(X,Y,U \mid Z)$, factorizes as :
\begin{eqnarray}
\label{factorisation}
P(X,Y,U \mid Z) = P(U) \; P(X \mid Z, U) \; P(Y \mid X, Z, U).
\end{eqnarray}
\noindent If we write $N(a,b)$
for the normal distribution  with
mean $a$ and variance $b$, our assumed model has form
\begin{eqnarray}
\label{full1}
P(U) &=& N(0,1),\\
\label{full2}
P(X \mid Z, U) &=& N(\omega_X + \alpha Z^{T} + \delta_X U, \sigma_X^2),\\
\label{full3}
P(Y \mid X, Z, U) &=& N(\omega_Y + \theta X + \beta Z^{T} + \delta_Y U,
\sigma_Y^2),
\end{eqnarray}
\noindent or, equivalently,
\begin{eqnarray}
\label{triangular1}
P(X \mid Z) &=& \omega_X + \alpha Z^{T} + A,\\
\label{triangular2}
P(Y \mid X, Z) &=& \omega_Y + \theta X + \beta Z^{T} + B.
\end{eqnarray}
\noindent  with $A \sim N(0, \delta_X^2+\sigma_X^2)$
$B \sim N(0, \delta_Y^2+\sigma_Y^2)$ and
$cov(A,B) = \delta_X \delta_Y$.
\newpar We exclude a possible effect of $(X,Z,U)$ on the variance
of $Y$. Equations (\ref{triangular1}--\ref{triangular2}) contain the
$2J+6$ parameters $(\omega_X,\omega_Y, \tau_X^2 \equiv \delta_X^2+\sigma_X^2,
\tau_Y^2 \equiv \delta_Y^2+\sigma_Y^2, \lambda \equiv \delta_X \delta_Y,
\theta,\alpha, \beta)$, which we shall refer to 
as the {\em structural parametrization}. Of inferential interest
among these  is parameter $\theta$ of
Equation (\ref{full3}), which can be interpreted in terms of change in $Y$
due to an intervention on $X$. Not all of the $2J+6$ structural parameters
are identified from the likelihood, for the following
reason. The conditional
expectation  $E(X \mid Z) = \omega_X + \alpha Z^{T}$
and the conditional variance $var(X \mid Z) = \tau_X^2$
provide $J+2$ conditions that make
structural parameters $\omega_X, \alpha, \tau_X$
identifiable from the likelihood.
The conditional expectation
$E(Y \mid X, Z) = \omega^{'}_Y + \theta^{'} X + \beta^{'} Z^{T}$
provides additional $J+2$ conditions: $\omega^{'}_Y =
\omega_Y - \omega_X \frac{\lambda}{\tau_X},
\theta^{'} = \theta + \frac{\lambda}{\tau_X}$ and 
$\beta^{'} = \beta - \alpha \frac{\lambda}{\tau_X}$.
A further condition is provided
by the conditional variance
$var(Y \mid X,Z) = \left(1-\left(\frac{\lambda}{
\tau_Y \tau_X}\right)^2\right)
(\tau_Y^2)$, for a total of $2J+5$ conditions -- insufficient
to estimate the full set of $2J+6$
structural parameters. The information contained in the data
identifies the structural parameters $\omega_X, \alpha, \tau_X$,
but fails to identify the remaining structural
parameters, $(\omega_Y, \tau_Y, \lambda,
\theta, \beta)$, including the
parameter of inferential interest, $\theta$. Full
parameter identifiability is achieved in the unlikely, and
therefore uninteresting, situation where
the values of the $J$ components of $\beta$ are
supplied by external knowledge.

\newpar We tackle the problem from
a Bayesian point of view, by designing a prior
distribution
that makes the posterior distribution proper. 
Let the symbol $D$ denote the data.
Then it is helpful to express the posterior
in the following product form:
\begin{eqnarray*}
\pi &\equiv& P(\omega, \tau, \lambda,
\theta,\alpha, \beta \mid D)\\
&=& P(\omega_X,\alpha,\tau_X \mid D) \;
P(\omega_Y,\theta, \tau_Y,\lambda  \mid \beta,
\omega_X,\alpha,\tau_X, D)
\;
P(\beta  \mid \omega_X,\alpha,\tau_X, D)
\end{eqnarray*}
\noindent Because $\beta$ belongs to the unidentifiable
subset of the model parameters, we have
$P(\beta  \mid \omega_X,\alpha,\tau_X, D)=
P(\beta  \mid \omega_X,\alpha,\tau_X)$,
which leads to
\begin{equation}
\label{final}
\pi = P(\omega_X,\alpha,\tau_X \mid D)
\; P(\omega_Y,\theta, \tau_Y, \lambda  \mid \beta, \omega_X, \alpha,\tau_X, D)
\;
P(\beta  \mid \omega_X,\alpha,\tau_X)
\end{equation}
\noindent The data allow us
to learn about parameters
$\omega_X,\alpha,\tau_X$ (first term of the above product) and,
conditional on $\beta, \omega_X,\alpha,\tau_X$, they
allow us to learn about the remaining parameters in the model
(second term of the product), but they provide no information
about $\beta$ and, in particular, about
the possible dependence between $\beta$ and $\alpha$. This is not a fatal flaw
if we provide $\beta$ with a suitable, scientifically plausible, prior.
One option consists of assuming $\beta =0$. Reasons why we 
repute this an untenable assumption have been previously discussed.
Another option is to impose inequalities based
on the assumption that, say, the direct component of the effect
of $Z_j$ on $Y$ is smaller in magnitude than the (indirect)
effect mediated by $X$. In this paper we adopt a
different approach, that requires us to
introduce a further condition that is sometimes referred
to as the {\em Instrument Strength Independent of Direct Effect}
(INSIDE) condition.

\begin{condition}[INSIDE]
\label{inside}
The genetic associations with the exposure
are independent of the direct effects of the genetic variants
on the outcome:
$\beta_j \indep \alpha_j$, for $j=1, \ldots , J$
\end{condition}

\newpar Under Condition \ref{inside}, we
shape the conditional prior
$P(\beta  \mid \omega_X,\alpha,\tau_X) = P(\beta)$
to express our prior belief that some of its components are
zero. This is discussed in the next section.

\subsection{The Prior}
\label{The Prior}

\newpar We shall now discuss the
prior specifications for
the parameters of model
(\ref{full1}--\ref{full3}).
Of special interest is the prior we impose
on the vector $\beta$. It is often assumed in the
MR literature that all the components of $\beta$ are zero. 
We replace this with the more realistic assumption that
{\em some} of the components of this vector are zero,
and incorporate this in our model by imposing on
$\beta$ the horseshoe shrinkage prior proposed by
Carvalho \cite{carvalho2010}. With this prior,
the components of $\beta$ will be shrunk
towards zero, but to different degrees
inferred from the data: large components will be only
moderately shrunk, while small components will be
heavily shrunk towards zero. The prior mechanism
may be informally described as follows. Those instruments
whose effects on $Y$ are irreconcilable with a no-pleiotropy
model will have large corresponding $\beta$ parameters. Our
prior will leave these parameters relatively unshrunk, so that
those instruments will have little impact on the estimated $\theta$.
By contrast, those instruments that are compatible with a
no-pleiotropy model and a low-variance outcome error
probability will have their corresponding $\beta$ parameters
heavily shrunk towards 0, so that the estimate of theta will
predominantly depend on the information provided by these
non-pleiotropic instruments.

\newpar We apply Carvalho's horseshoe 
prior to $\beta$ by writing:
\begin{eqnarray}
\label{scale mixture 1}
p(\beta_j \mid \phi_j) &=& N(0, \phi_j^2),\\
\label{scale mixture 2}
p(\phi_j \mid \gamma) &=& Cauchy^{+}(0,\gamma),\\
\label{scale mixture 3}
p(\gamma) &=& Cauchy^{+}(0,1
),
\end{eqnarray}
\noindent for $j=1, \ldots , J$, where
$Cauchy^{+}(0,a)$ denotes the
half-Cauchy distribution
on the positive reals with scale
parameter $a$.
Crucially, in the above prior, 
each $\beta_j$ is mixed over
its own unknown $\phi_j$. The parameters
$\phi_j$s control the {\em local} degrees
of shrinking, and are all independently
drawn from a half-Cauchy prior with a common,
unknown global scale parameter $\gamma$, which 
controls the {\em global} degree of shrinking.
By virtue of (\ref{scale mixture 1}),
small values of $\phi_j$ cause
$\beta_j$ to shrink towards zero,
whereas large values will
prevent the estimate of $\beta_j$
from shrinking.
Importantly, the horseshoe
prior is free from user-chosen
hyperparameters.

\newpar The shrinkage for instrument $j$
is measured by
parameter $\kappa_j = 1/(1+ \phi_j^2)$, called
the shrinkage weight, with
$\kappa_j=0$ (resp., $\kappa_j=1$)
indicating absence of (resp., near-total) shrinking.
The shrinkage weights $\kappa_j$ are
inferred from the data. Equations (\ref{scale mixture 1}--
\ref{scale mixture 3}),
with $\gamma=1$, yield a horseshoe-shaped $Beta(.5,.5)$
prior for $\kappa_j$, peaking at $\kappa_j=0$
and $\kappa_j=1$. The ability of our model to 
discriminate between pleiotropic and non-pleiotropic instruments
corresponds to the tendency of parameter $\kappa_j$ to
be, on average, lower for the non-pleiotropic than for the pleiotropic
instruments. In Section \ref{Simulation Experiment}
we assess separation by a simulation experiment.

\newpar We are now going to discuss the
prior specifications for the remaining model parameters.
In our analyses, we have taken parameters $\omega_X$
and $\omega_Y$ to follow a priori independent uniform
distributions. However, since these parameters are related
to the global means of two observed variables, one may be
able to shape informative priors for $\omega_X$ and
$\omega_Y$ on the basis of external information. 
In our analyses we took each $\alpha_j$, for $j=1, \ldots , J$,
to be independently drawn from a normal $N(\mu_{\alpha},
\sigma_{\alpha}^2)$ population prior, with hyperparameters
$\mu_{\alpha}$ and $\sigma_{\alpha}$ subject 
to uniform priors. However, under our $Z  \indep U$ assumption
one will often be able to shape an informative prior
for the $\alpha$ parameters (for example a prior that imposes
on these parameters specific signs) on the basis of external $(Z,X)$ data.
Parameters  $\sigma_X$ and $\sigma_Y$ are only partially
identifiable. We assigned these
two parameters uniform positive prior distributions,
which did not cause mixing problems
in our Markov chain Monte Carlo exploration
of the model posterior distribution.  We completed our prior specifications
by taking the parameters $\eta_X$, $\eta_Y$ and $\theta$ to follow
uniform independent priors.

\newpar There are situations where some instruments can be safely
assumed to exert no pleiotropic effect. These instruments can be used
to gather prior information about the value of  $\tau_X$. By combining
this information with the constraint $\tau_X^2 = \eta_X^2+\sigma_X^2$,
one may attempt to derive inequalities involving
$\eta_X$ and $\sigma_X$ and use them to shape an informative
joint prior for these two parameters.

\begin{figure}
\begin{center}
\hspace*{-0cm}\includegraphics[scale=0.6]{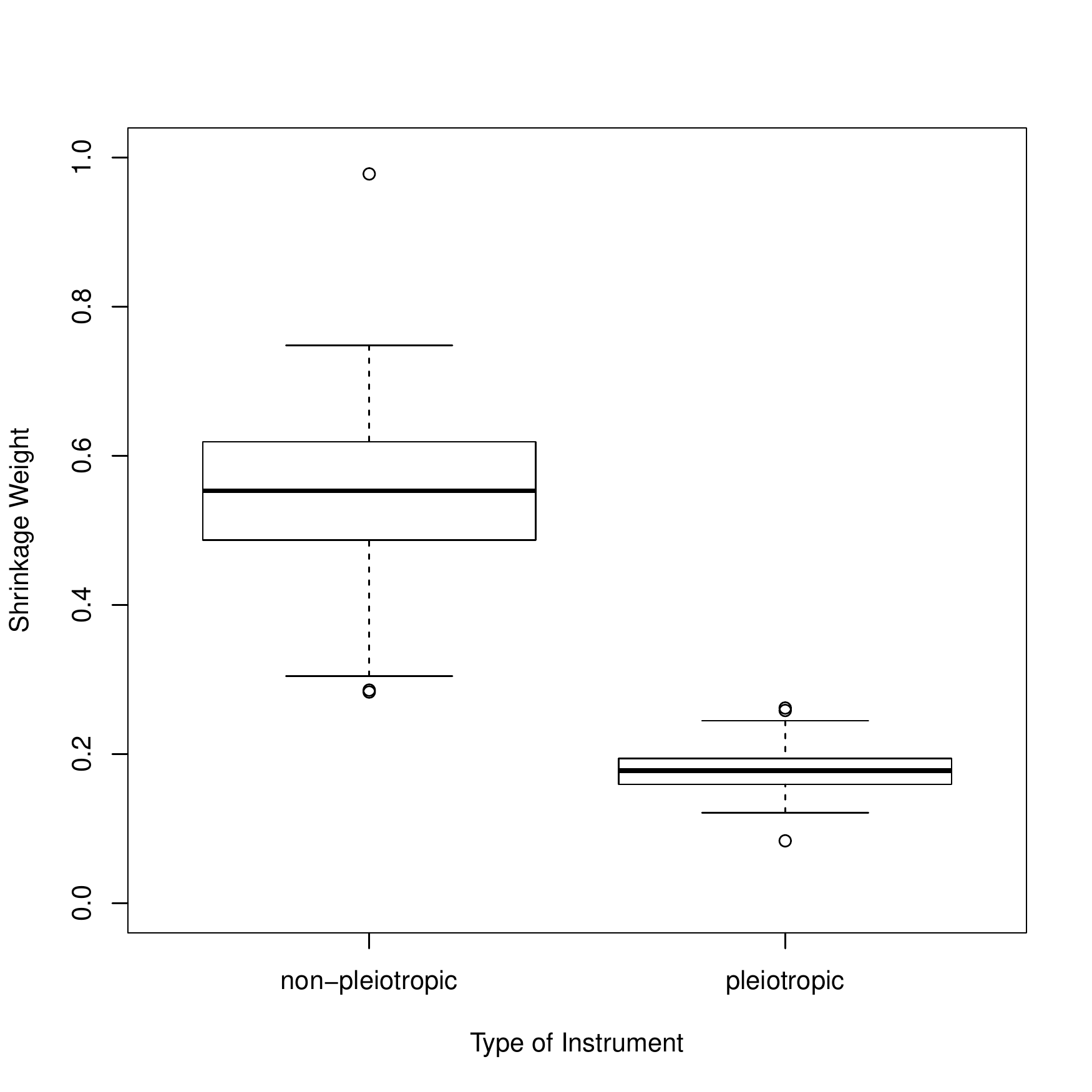}
\caption{\small We have selected at random one of the datasets
generated in the simulation experiment of Section \ref{Simulation Experiment}, and compared
the distribution of the shrinkage parameter $\kappa_j$ over the set
of the non-pleiotropic instruments (box on the left) with that over
the set of pleiotropic instruments (box on the right).
As expected, the components of $\beta$ associated with
the pleiotropic instruments have values of $\kappa_j$ close to zero,
which indicates that they tend to be left unshrunk,
whereas the components of $\beta$ associated with
the non-pleiotropic instruments have values of $\kappa_j$ closer to 1,
which indicates that they are heavily shrunk.}
\end{center}
\end{figure}

\section{Simulation Experiment}
\label{Simulation Experiment}

A simulation experiment was set up to comparatively assess the
performance of the proposed method and of the weighted median
estimator (WME) \cite{Han2008} \cite{Bowden2016}. 
Throughout the experiment,
we set the number of instruments, $J$, to be equal to 20, with
the instrumental variables 
representing "allele doses" (0, 1, 2), and we took the effect
of each instrument on both $X$ and $Y$ to
be linear in the allele dose. We imposed a pleiotropic effect
on half of the instruments, the remaining instruments
being treated as pleiotropic.
We considered six simulation
scenarios differing by type of pleiotropy
(balanced, positive
or negative) and sample size
(100, 520). For each scenario, we
simulated 1000 replicate datasets by setting $\theta=0$
-- the null hypothesis --,
and further 1000 datasets by setting
$\theta= 0.35$ -- the alternative hypothesis. Each new simulation started
with the generation of 
a configuration of values of
$Z_1, \ldots, Z_{20}$ for each
hypothetical individual, these values being
drawn from $J$
independent trinomial distributions that
mimicked the joint distribution of real SNP loci.
Conditional on such values, the simulation
proceeded with the generation of a configuration
of values of $(X,Y,U)$ 
for each hypothetical individual,
on the basis of Equations (\ref{full1}--\ref{full3})
with
$\delta_X \sim N(-.05, .0025),
\delta_Y \sim N(-.1, .0025),
\omega_Y \sim N(-3.7, .04),
\omega_X = 3.3,
\sigma_Y = 0.1,
\sigma_X = 0.1$ and
$\alpha_j   \sim N(.034, .0031),
\beta_j \sim N(.012 \xi, .0025)$ , for $j=1, \ldots, 10$, where
$\xi$ indicates balanced ($\xi=0$), negative ($\xi=-1$) and
positive ($\xi=1$) pleiotropy, and $\beta_j =0$ , for $j=11, \ldots, 20$.

\newpar Each simulated dataset
was analyzed by using both the WME and
the proposed method.
In the latter case, inference was based on
Markov chain Monte Carlo (MCMC)
samples from the posterior distribution defined
by Equations (\ref{full1}--\ref{full3}) 
jointly with the
prior specifications of Section \ref{The Prior}. These samples
were generated by using the Hamiltonian dynamics
MCMC techniques
\cite{Metropolis1953} \cite{Neal2011}
offered by the program {\tt Stan}
\cite{STAN2014a}\cite{STAN2014b}.
     Initial values for the Markov chains were generated
automatically in STAN based on approximate posterior mean
estimates obtained via variational inference techniques
\cite{Wainwright2008}.
No Markov chain mixing problems were encountered.
The WME analysis of each dataset produced
a point estimate and a corresponding
95 percent confidence interval for $\theta$. Analysis by
our method produced a
posterior mean and a 95 percent Bayesian credible interval 
for $\theta$.

\newpar The two methods were comparatively assessed in terms
of {\em (i)} coverage under the null, {\em (ii)} coverage under the alternative,
{\em (iii)} power, {\em (iv)} bias under the null
and {\em (v)} bias under the alternative. The estimated
coverage under the null (resp., alternative) was
the proportion of simulations performed under
the null (resp., alternative) where the credible/confidence 
interval for $\theta$ contained the value $0$ (resp., $0.35$). 
Power was estimated as the proportion of simulations under
the alternative  hypothesis where the credible/confidence interval 
for $\theta$ was contained by the positive real axis.
Bias was estimated as the average  signed difference between
the point estimate for $\theta$ and the corresponding true
value of the parameter.
A summary of the simulation results is 
given in Table 1. The proposed method appears
to outperform the WME in terms of coverage under the null, especially
in the presence of a larger sample size. It appears superior also
in terms of coverage under the alternative,
when the sample is small. The two
methods appear to offer roughly the same power. As expected, in
both methods, power appears to depends on sample size. The differences
in bias between the two methods are minimal. Under the proposed method,
coverage under the null appears to be
slightly more robust to
positive pleiotropy than with the WME, when the sample is small.

\begin{table}[ht!]
\centering
\begin{tabular}{@{}ccccccccccccccc@{}}
&&&\phantom{ab}& \multicolumn{5}{c}{Our method} &
\phantom{ab}& \multicolumn{5}{c}{Weighted median}\\
&&&\phantom{ab}& \multicolumn{5}{c}{} &
\phantom{ab}& \multicolumn{5}{c}{estimator}\\
\cmidrule{5-9} \cmidrule{11-15}\\
\\
\\
\\
\\
\\
\\
\\
\\
\begin{rotate}{90}Scenario \end{rotate}&
\begin{rotate}{90}Pleiotropy\end{rotate}&
\begin{rotate}{90}Sample Size \end{rotate}&&
\begin{rotate}{90}Coverage under Null \end{rotate}&
\begin{rotate}{90}Coverage under Alternative \end{rotate}&
\begin{rotate}{90}Power\end{rotate}&
\begin{rotate}{90}Bias under Null\end{rotate}&
\begin{rotate}{90}Bias under Alternative \end{rotate}&&
\begin{rotate}{90}Coverage under Null \end{rotate}&
\begin{rotate}{90}Coverage under Alternative \end{rotate}&
\begin{rotate}{90}Power\end{rotate}&
\begin{rotate}{90}Bias under Null\end{rotate}&
\begin{rotate}{90}Bias under Alternative \end{rotate}\\
\midrule
1&0&520&&.92&.9&.89&.001&.025&&.79&.79&.875&-.002&.007\\
2&$-$&520&&.89&.9&.86&-.028&.008&&.77& .80& .85&-.043&-.009\\
3&$+$&520&&.90&.86& .89&.029&-.017&&.79&.77&.91&.038&.015\\
4&0&100&&.90&.92&.57&.004&.3&&.88&.93&.53&.03&.055\\
5&$-$&100&&.92&.91& .54&-.32&-.024&&.89& .92& .55&.009&-.02\\
6&$+$&100&&.89&.92&.57&.029&.06&&.81&.84&.65&.07&.06\\
\hline
\end{tabular}
\caption{\small
Results
of a comparative assessment of the
proposed and of the median estimator
methods for causal effect estimation.
Table rows correspond to six different
simulation scenarios characterized
by the presence of balanced (0), positive ($+$)
or negative ($-$) pleiotropy and by the sample
size (520 vs 100). Throughout the simulation,
the number of instruments was
kept equal to 20. Method performance
was separately assessed under each of the 6 scenarios
on the basis of 2000 simulated datasets, in terms
of coverage under the null and under the alternative,
power and bias. See main text for further details.}
\end{table}

\newpar In Figure 3 we have selected at random one of the datasets
generated in this simulation experiment, and compared
the distribution of the shrinkage parameter $\kappa_j$ over the set
of the non-pleiotropic instruments (box on the left) with that over
the set of pleiotropic instruments (box on the right).
As expected, the components of $\beta$ associated with
the pleiotropic instruments have values of $\kappa_j$ close to zero,
which indicates that they tend to be left unshrunk,
whereas the components of $\beta$ associated with
the non-pleiotropic instruments have values of $\kappa_j$ closer to 1,
which indicates that they are heavily shrunk.

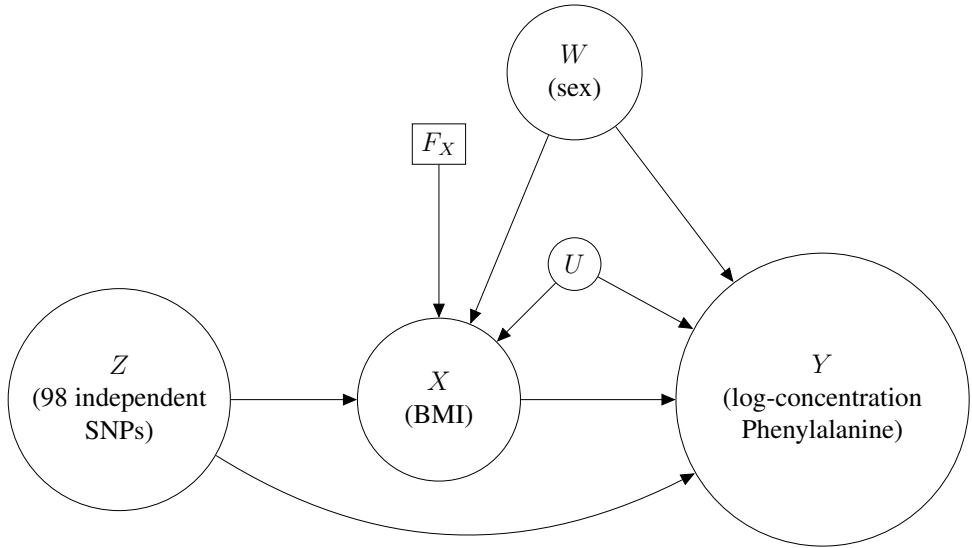
\begin{figure}[!h]
\begin{center}
\scalebox{.85}{
\centering
\begin{tikzpicture}[align=center,node distance=7cm]
\node (X) [draw=black,fill=none,circle] {\begin{minipage}{2cm}\begin{center}$X$\\(BMI)\end{center}\end{minipage}};
\node (U) [above right of=X,draw=black,fill=none,circle,
node distance=3cm] {$U$};
\node (W) [above of=U,draw=black,fill=none,circle,
node distance=3cm] {\begin{minipage}{1.5cm}\begin{center}$W$\\(sex)\end{center}\end{minipage}};
\node (Y) [right of=X,draw=black,fill=none,
circle,node distance=6cm] {\begin{minipage}{4cm}\begin{center}
$Y$\\(log-concentration Phenylalanine)
\end{center}\end{minipage}};
\node (Z) [left of=X,draw=black,fill=none,circle,,node distance=5cm] {\begin{minipage}{2.8cm}\begin{center}$Z$\\(98 independent SNPs)\end{center}\end{minipage}};
\node (SX) [above of=X,draw=black,fill=none,rectangle,node distance=4cm]{$F_X$};
\draw[-triangle 45] (X) to node
[auto] {} (Y);
\draw[-triangle 45] (U) to node {} (X);
\draw[-triangle 45] (U) to node {} (Y);
\draw[-triangle 45] (Z) to node {} (X);
\draw[-triangle 45] (SX) to node {} (X);
\draw[-triangle 45] (W) to node {} (X);
\draw[-triangle 45] (W) to node {} (Y);
\draw[-triangle 45, bend right] (Z) to node {} (Y);
\end{tikzpicture}\label{Figure 1}
}
\caption{\small Structure of the illustrative problem.}
\end{center}
\end{figure}

\begin{figure}
\hspace*{-0cm}\includegraphics[scale=0.8]{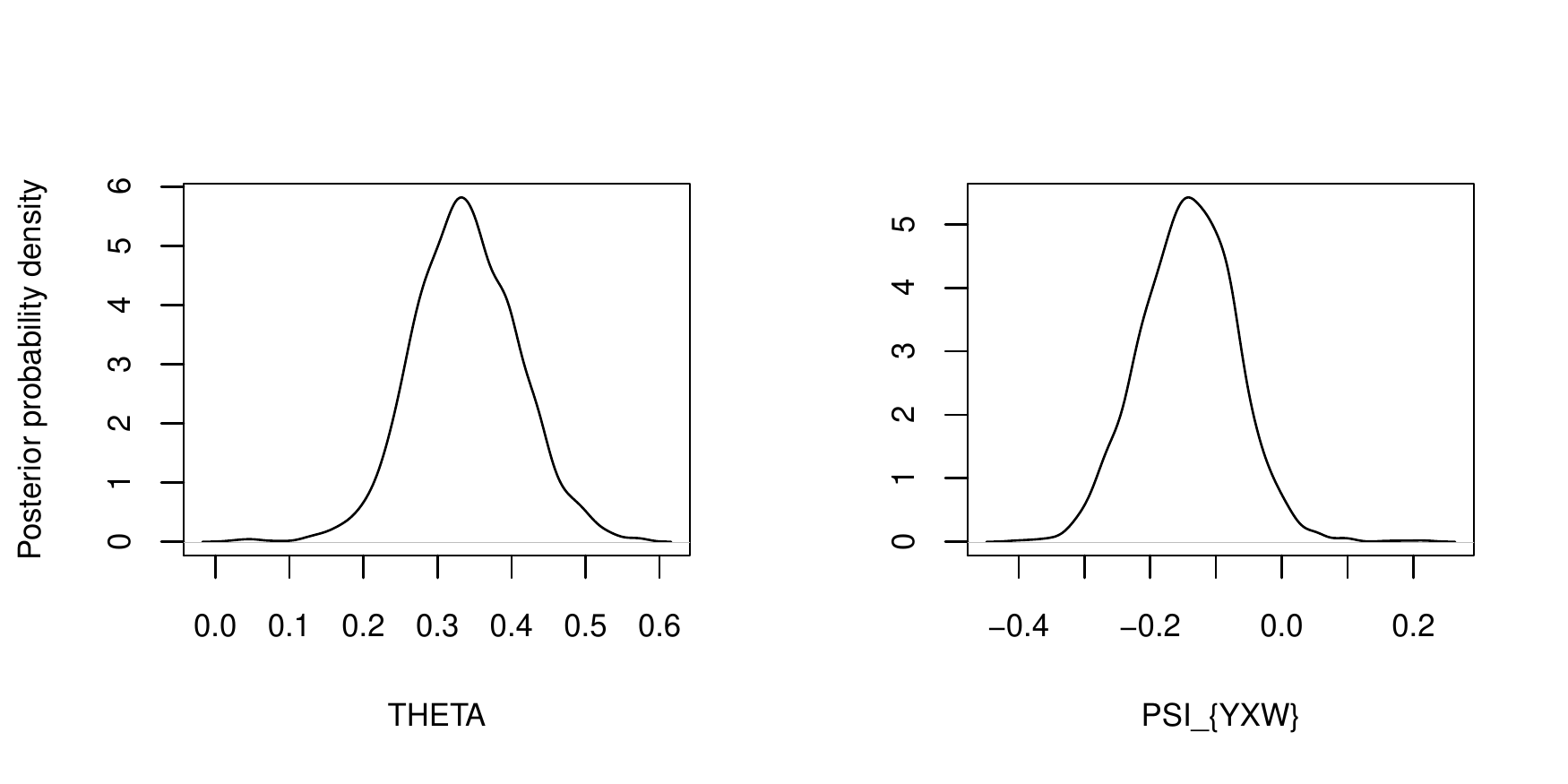}
\caption{\small Marginal posterior densities for parameters $\theta$
and $\psi_{YXW}$ of Equations (\ref{fullinteraction}), based on the DILGOM data of
Section \ref{Sex-Dependent Causal Effect of Body Mass on Phenylalanine}}
\end{figure}

\begin{figure}
\hspace*{-0cm}\includegraphics[scale=0.7]{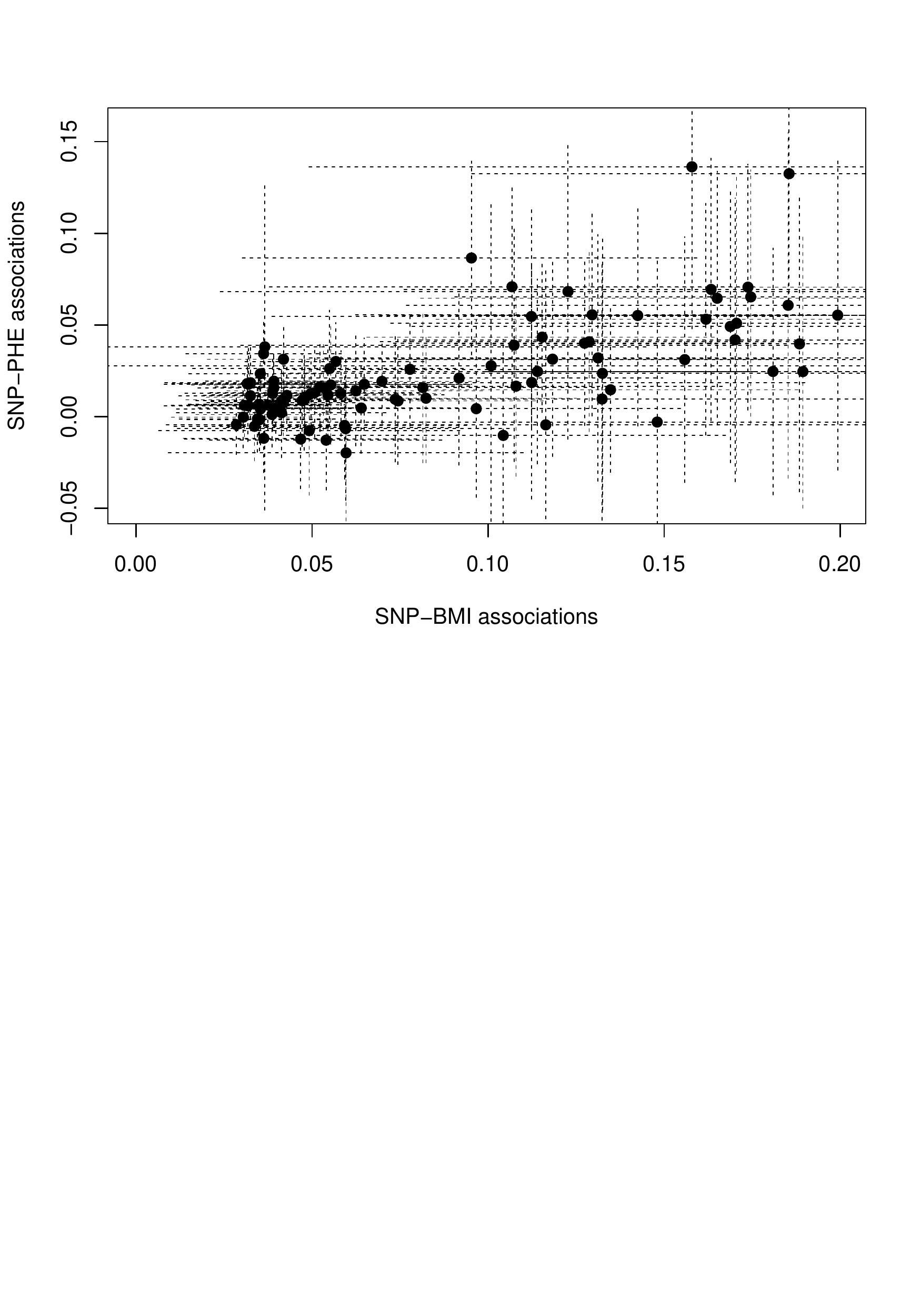}
\vspace{-9cm}
\caption{\small In this plot,
the black dots
correspond to the instrumental SNPs
in the analysis of Section
\ref{Sex-Dependent Causal Effect of Body Mass on Phenylalanine}, the
horizontal coordinates to the SNP's coefficient
in  the exposure regression (least-squares regression
of BMI on that SNP),
and their vertical coordinates to the SNP's coefficient
in  the response regression. The 95 percent credible
intervals for the coefficients are represented as dashed segments.}
\end{figure}

\section{Sex-Dependent Causal Effect of Body Mass on Phenylalanine}
\label{Sex-Dependent Causal Effect of Body Mass on Phenylalanine}

\noindent A high body mass is associated with an increased
risk of several chronic diseases. A better understanding of the 
underlying biology requires that we identify the metabolic
mediators of this deleterious effect \cite{moore2014}.
Within this perspective, Ho \etal \cite{ho2016}
have analyzed data from 2383 Framingham offspring cohort
participants, and tested the association between
the body mass index (BMI) and
more than of two hundred
cardiometabolic traits and metabolites. As many as sixty
metabolites were found to be significantly (P $<$ 0.00023)
associated with BMI. The next step is now to assess
whether these associations are causal.
We advocate a MR approach to the problem, where BMI
acts as exposure and the individual metabolites,
in turn, act as responses. SNPs associated with BMI
are used as instruments to assess 
whether BMI is a causal influence on the metabolite
of interest. 

\newpar The idea can be effectively implemented by
using our proposed method. To
illustrate this, we shall now apply our method to the assessment of the putative causal effect of BMI on one of the
metabolites highlighted by the Framingham study: the aromatic amino acid phenylalanine. 
This analysis we have carried out
on the basis of data from 520 unrelated individuals, aged 25--74 years, sampled from a population-based Finnish cohort --
the Dietary, Lifestyle and Genetic determinants of Obesity and Metabolic Syndrome (DILGOM)
study\cite{inouye2010}. The DILGOM data contain  information about the individuals' serum metabonomes,
combined with the genome-wide profiles of genetic and transcriptional variation from blood leukocytes of the same individuals. 
The data also contain individual-level measures of BMI, age and sex. 

\newpar The first step of our analysis consisted of selecting a subset of 98 DILGOM-genotyped
SNPs, based  on a $p \leq 1 \times 10^{-5}$ significance
threshold for the SNP's  association with BMI  (more precisely, for the Wald test
statistic for the regression of BMI on the SNP),
and on a maximum between-SNP  linkage disequilibrium score of $\leq 0.05$. We let
the genotypes of the 98 selected SNPs act as instruments
in our analysis, $(Z_1, \ldots , Z_{98})$, to be treated as continuous variables, with
values $(0,1,2)$ corresponding to
the number of minor alleles found at the SNP. We let the variable BMI act as exposure, $X$,
and the concentration of
phenylalanine, expressed on a log scale, act as outcome, $Y$. We incorporate in the analysis
the variable sex, denoted as $W$,
taking value $0$ for female and $1$ for male. The resulting problem
structure is depicted in Figure 4.

\newpar The model specification (\ref{full1}-\ref{full3}) was
modified to acknowledge the possible
interaction between the effects of sex ($W$) and 
BMI ($X$) on phenylalanine concentration ($Y$). 
Such an elaboration was motivated by evidence \cite{kaplan2014}
that sex and BMI interact in their effects on 
coronary artery disease and other clinical outcomes.
The elaboration
was straightforwardly implemented, as is often
the case within a Bayesian approach that keeps model specification and
inference calculations independent of each other. All we had to
do was to modify Equations (\ref{full1}-\ref{full3}) into: 
\begin{eqnarray}
P(U) &=& N(0,1),\nonumber\\
\label{fullinteraction}
P(X \mid Z, U) &=& N(\omega_X + \alpha Z^{T} + \psi_{XW} W + \delta_X U, \sigma_X^2),\\
P(Y \mid X, Z, U) &=& N(\omega_Y + \theta (X+ \psi_{YXW} W) + \beta Z^{T} + \psi_{YW} W+ \delta_Y U,
\sigma_Y^2)\nonumber,
\end{eqnarray}
\noindent where the $BMI \times SEX$ interaction is 
represented by parameter $\psi_{YXW}$, so that the
causal effect of BMI on
log-concentration of phenylalanine is represented by
$\theta$ in the female stratum, and by
$\theta^{'} \equiv \theta+\psi_{YXW}$ in the males. 

\newpar We completed the model by specifying the priors
as described in Section \ref{The Prior}. We sampled the model posterior distribution
by running 10000 iterations of a (Hamiltonian dynamics) Markov chain.
The last 5000 samples were used to approximate the posterior
marginal distributions and to obtain posterior
means and credible intervals for the parameters of interest.
Figure 5 shows the obtained marginal
posterior distributions for parameters $\theta$
and $\psi_{YXW}$ of Equations (\ref{fullinteraction}).
The parameter $\theta^{'}$ was included in the model as a derived quantity,
so as to obtain samples from its marginal posterior distribution
and to estimate its posterior mean and credible interval.
 The estimates for the parameters of interest are reported in Table 1.
The 95 percent credible interval for $\psi_{YXW}$ is contained by the positive real axis,
which represents strong evidence of an interaction between the causal effects of
BMI and sex on the concentration of phenylalanine. The causal
effect of BMI on the log-concentration of phenylalanine was estimated
to be 0.34 (95 percent credible
interval 0.21 to 0.47) in the females,
and 0.2 (95 percent credible interval 0.098 to 0.3)
in the males.

\newpar It would have been
possible to examine the interaction by
running separate analyses within
the male and female strata, although and the cost of
a loss of statistical power. Figure 6 provides visual evidence
of the causal effect of the BMI on phenylalanine. In this figure,
each instrumental SNP appears as a black dot with the $x$ and
$y$ coordinates corresponding to the coefficients of the SNP in
the exposure and outcome regressions, respectively, and the
corresponding 95\% confidence intervals are represented by
dashed segments. The emerging linear relationship can be
interpreted to suggest that hypothetical perturbations of $X$
would result in corresponding, proportional, perturbations in $Y$.

\newpar By repeating the analysis on a large set of metabolites, we
may aim to a classification of obesity based on the values of the
molecular mediators that are responsible for its deleterious effects,
for purposes of personalised medicine and drug target discovery.

\begin{table}[ht!]
\centering
\begin{tabular}{@{}cccc@{}}\\
&
Posterior&
Standard&
\\
Parameter&
Mean&
Deviation&
95 Percent Credible Interval\\
\midrule
$\theta$& 0.34 & 0.068&(0.21,0.47)\\ 
$\psi_{YXW}$&-0.14 &0.072&(-0.28,-0.0062)\\
$\theta^{'}$& 0.2&0.051&(0.098,0.3)\\ 
\hline
\end{tabular}
\caption{\small
Point and interval estimates for some parameters
of model (\ref{fullinteraction}), based on the DILGOM
dataset of Section \ref{Sex-Dependent Causal Effect of Body Mass on Phenylalanine}.
The causal effect of BMI on the
log-concentration of phenylalanine is represented by parameter
$\theta$ in the female stratum, and by
$\theta^{'} \equiv \theta+\psi_{YXW}$ in the males.}
\end{table}

\section{Discussion}
\label{Discussion}

Kang and colleagues \cite{Kang2016} use the
term "invalid" (resp., "valid") to denote an instrument
that violates (resp., obeys) the Exclusion Restriction condition.
They propose a LASSO-type procedure to identify the
valid instruments from within a set of candidate instrumental
variables.  The idea, further elaborated by Windmeijer and
colleagues \cite{Windmeijer2016}, is to obtain a sparse estimate
of the vector representing the pleiotropic effects by imposing on
it an $l_1$ penalty. The null elements of the estimated vector
should then correspond to the valid instruments. In our framework,
we may construct a Bayesian analogue of Windmeijer's approach
by replacing the horseshoe prior on $\beta$ with a double-exponential
prior. But this will make a big difference to the posterior means
when $\beta$ is sparse. This is because our (horseshoe) prior
presents superior tail robustness to the large signals introduced in
the $\beta$ vector by the pleiotropic effects, and possesses the ability
to shrink the components of $\beta$ near zero much more forcefully
than those far from zero, thanks to the combined local and global shrinking.
Our work differs from Windmeijer's also from the point of view of the
method justification, which in the present paper is based on Dawid's
decision-theoretic causal inference framework \cite{Dawid2015}, rather
than on the notion of potential outcome.

\newpar A number of method issues await investigation.
Of foremost importance is to consider the method behaviour
in the presence of  collinearity of the
instruments in the outcome regression. Our prior on $\beta$ should
supply enough tail probability to produce a posterior distribution 
which preserves the pattern of correlation between the components
of $\beta$.

\section*{Acknowledgments}

The first two authors acknowledge partial support from the European
Union's Seventh Framework Programme FP7-Health- 2012-INNOVATION, under grant agreement number 305280 (MIMOmics).
The DILGOM data resource exploited in Section
\ref{Sex-Dependent Causal Effect of Body Mass on Phenylalanine} has been 
funded by the Sigrid Juselius and Yrj\~o Jahnsson Foundations
and by the Finnish Academy grants no. 255935 and 269517. 
Our analysis of these data has benefitted from discussions
with Drs. Xiaoguang Xu and
Susana Conde.

\bibliographystyle{plain}
\bibliography{mendelian}

\begin{thebibliography}{10}

\bibitem{Bowden2015}
J~Bowden, G~Davey~Smith, and S~Burgess.
\newblock Mendelian {R}andomization {W}ith {I}nvalid {I}nstruments: {E}ffect
  {E}stimation and {B}ias {D}etection {T}hrough {E}gger {R}egression.
\newblock {\em International Journal of Epidemiology}, 44(2):512–525, 2015.

\bibitem{Bowden2016}
Jack Bowden, George Davey~Smith, Philip~C Haycock, and Stephen Burgess.
\newblock Consistent estimation in mendelian randomization with some invalid
  instruments using a weighted median estimator.
\newblock {\em Genetic Epidemiology}, 40:304--314, 2016.

\bibitem{Burgess2012}
S~Burgess, A~Butterworth, A~Malarstig, and S.~Thompson.
\newblock Use of {M}endelian {R}andomisation to {A}ssess {P}otential {B}enefit
  of {C}linical {I}ntervention.
\newblock {\em BMJ}, 345:e7325, 2012.

\bibitem{burgess2016sensitivity}
Stephen Burgess, Jack Bowden, Tove Fall, Erik Ingelsson, and Simon~G Thompson.
\newblock Sensitivity analyses for robust causal inference from {M}endelian
  randomization analyses with multiple genetic variants.
\newblock {\em Epidemiology (in press)}, 2016.

\bibitem{Burgess2015}
Stephen Burgess and Simon~G. Thompson.
\newblock {\em {M}endelian {R}andomization: Methods for Using Genetic Variants
  in Causal Estimation.}
\newblock Chapman and Hall, 2015.

\bibitem{Han2008}
Han C.
\newblock Detecting invalid instruments using l1-gmm.
\newblock {\em Econ Lett}, 101:285--287, 2008.

\bibitem{carvalho2010}
Carlos~M Carvalho, Nicholas~G Polson, and James~G Scott.
\newblock The horseshoe estimator for sparse signals.
\newblock {\em Biometrika}, page asq017, 2010.

\bibitem{Davey2003}
G~Davey~Smith and S.~Ebrahim.
\newblock Mendelian {R}andomization: {C}an {G}enetic {E}pidemiology
  {C}ontribute to {U}nderstanding {E}nvironmental {D}eterminants of {D}isease?
\newblock {\em International Journal of Epidemiology}, 32:1–22, 2003.

\bibitem{Davey2014}
G~Davey~Smith and G.~Hemani.
\newblock Mendelian {R}andomization: {G}enetic {A}nchors for {C}ausal
  {I}nference in {E}pidemiological {S}tudies.
\newblock {\em Hum Mol Genet}, 23:89--98, 2014.

\bibitem{Dawid1979}
A.~P. Dawid.
\newblock Conditional independence in statistical theory (with {D}iscussion).
\newblock {\em J. Roy. Statist. Soc. B}, 41:1--31, 1979.

\bibitem{Dawid2002}
A.~Philip Dawid.
\newblock Influence diagrams for causal modelling and inference.
\newblock {\em International Statistical Review}, 70:161--189, 2002.
\newblock Corrigenda, {\em ibid\/}., 437.

\bibitem{Dawid2015}
A.~Philip Dawid.
\newblock Statistical causality from a decision-theoretic perspective.
\newblock {\em Ann. Rev. Statist. Appl.}, 2:273--303, 2015.

\bibitem{Didelez2007}
Vanessa Didelez and Nuala~A. Sheehan.
\newblock Mendelian randomisation as an instrumental variable approach to
  causal inference.
\newblock {\em Statistical Methods in Medical Research}, 16:309--330, 2007.

\bibitem{Didelez2007b}
Vanessa Didelez and Nuala~A. Sheehan.
\newblock Mendelian randomisation: Why epidemiology needs a formal language for
  causality.
\newblock In Federica Russo and Jon Williamson, editors, {\em Causality and
  Probability in the Sciences}, volume~5 of {\em Texts In Philosophy Series},
  pages 263--292. College Publications, London, 2007.

\bibitem{geiger90}
D.~Geiger, T.~Verma, and J.~Pearl.
\newblock Identifying independence in {B}ayesian networks.
\newblock {\em Networks}, 20(5):507--534, 1990.

\bibitem{ho2016}
JE~Ho, MG~Larson, A~Ghorbani, S~Cheng, MH~Chen, M~Keyes, EP~Rhee, CB~Clish,
  RS~Vasan, RE~Gerszten, and Wang TJ.
\newblock Metabolomic {P}rofiles of {B}ody {M}ass {I}ndex in the {F}ramingham
  {H}eart {S}tudy {R}eveal {D}istinct {C}ardiometabolic {P}henotypes.
\newblock {\em PloS one}, 11, 2016.

\bibitem{inouye2010}
M.~Inouye, J.~Kettunen, P.~Soininen, K.~Silander, S.~Ripatti, L.S. Kumpula,
  E.~Hamalainen, P.~Jousilahti, A.J. Kangas, S.~Mannisto, M.J. Savolainen,
  A.~Jula, J.~Leiviska, A.~Palotie, V.~Salomaa, M.~Perola, M.~AlaKorpela, and
  L.~Peltonen.
\newblock Metabonomic, transcriptomic, and genomic variation of a population
  cohort.
\newblock {\em Molecular Systems Biology}, 6:502–513, 2010.

\bibitem{Jones2012}
E.~M. Jones, J.~R. Thompson, V.~Didelez, and N.~A. Sheehan.
\newblock On the choice of parameterisation and priors for the bayesian
  analyses of mendelian randomisation studies.
\newblock {\em Statistics in Medicine}, 31(14):1483--1501, 2012.

\bibitem{Kang2016}
Hyunseung Kang, Anru Zhang, T.~Tony Cai, and Dylan~S. Small.
\newblock Instrumental variables estimation with some invalid instruments and
  its application to mendelian randomization.
\newblock {\em Journal of the American Statistical Association},
  111(513):132--144, 2016.

\bibitem{kaplan2014}
Robert~C. Kaplan, M.~Larissa Avilés-Santa, Christina~M. Parrinello, David~B.
  Hanna, Molly Jung, Sheila~F. Castañeda, Arlene~L. Hankinson, Carmen~R.
  Isasi, Orit Birnbaum-Weitzman, Ryung~S. Kim, Martha~L. Daviglus, Gregory~A.
  Talavera, Neil Schneiderman, and Jianwen Cai.
\newblock Body mass index, sex, and cardiovascular disease risk factors among
  hispanic/latino adults: Hispanic community health study/study of latinos.
\newblock {\em Journal of the American Heart Association}, 3(4), 2014.

\bibitem{lauritzen90}
S.~L. Lauritzen, A.~P. Dawid, B.~N. Larsen, and H.~G. Leimer.
\newblock Independence properties of directed {M}arkov fields.
\newblock {\em Networks}, 20(5):491--505, 1990.

\bibitem{Metropolis1953}
N.~Metropolis, A.~Rosenbluth, M.~Rosenbluth, M.~Teller, and E.~Teller.
\newblock Equations of state calculations by fast computing machines.
\newblock {\em Journal of Chemical Physics}, 21:1087--1092, 1953.

\bibitem{moore2014}
Steven~C. Moore, Charles~E. Matthews, Joshua~N. Sampson, Rachael~Z.
  Stolzenberg-Solomon, Wei Zheng, Qiuyin Cai, Yu~Ting Tan, Wong-Ho Chow,
  Bu-Tian Ji, Da~Ke Liu, Qian Xiao, Simina~M. Boca, Michael~F. Leitzmann, Gong
  Yang, Yong~Bing Xiang, Rashmi Sinha, Xiao~Ou Shu, and Amanda~J Cross.
\newblock Human metabolic correlates of body mass index.
\newblock {\em Metabolomics}, 10:0, 2013.

\bibitem{Neal2011}
R.~Neal.
\newblock {MCMC} using hamiltonian dynamics.
\newblock In Gelman A. Jones G.~L. Brooks, S. and X.L. Meng, editors, {\em
  Handbook of Markov Chain Monte Carlo}, pages 116--162. Chapman and Hall/CRC,
  2011.

\bibitem{STAN2014b}
Stan~Development Team.
\newblock {\em RStan, version 2.2.}
\newblock http://mc-stan.org/rstan.html, 2014.

\bibitem{STAN2014a}
Stan~Development Team.
\newblock {\em Stan: A C++ library for probability and sampling, version 2.2}.
\newblock http://mc-stan.org/, 2014.

\bibitem{Wainwright2008}
M.~J. Wainwright and M.~I. Jordan.
\newblock Graphical models, exponential families, and variational inference.
\newblock {\em Foundations and Trends in Machine Learning}, 1(1-2):1--305,
  2008.

\bibitem{Windmeijer2016}
Frank Windmeijer, Helmut Farbmacher, Neil Davies, and George Davey~Smith.
\newblock On the use of the {Lasso} for {I}nstrumental {V}ariables {E}stimation
  with {S}ome {I}nvalid {I}nstruments.
\newblock Bristol economics discussion papers, Department of Economics,
  University of Bristol, UK, 2016.

\end{thebibliography}
\end{document}